\documentclass[a4paper,twoside,11pt]{article}
\usepackage{a4wide} 
\usepackage[utf8]{inputenc}
\usepackage[T1]{fontenc}
\usepackage{lmodern}

\usepackage{graphicx} 
\usepackage{amsmath,amssymb,amsthm}
\usepackage{enumerate}
\usepackage{array}
\usepackage{multirow}

\usepackage{caption}
\usepackage{subcaption}

\usepackage[usenames,dvipsnames]{xcolor}
\definecolor{celestialblue}{rgb}{0.29, 0.59, 0.82}

\usepackage{url}
\usepackage[noadjust]{cite}
\usepackage{booktabs}

\usepackage{helvet}

\usepackage{tikz}
\usetikzlibrary{calc}
\usetikzlibrary{arrows}
\usetikzlibrary{decorations.pathreplacing}  
\usepgflibrary{shapes.geometric}

\usepackage{titling}
\usepackage{fancyhdr}
\usepackage[scale=1.0]{ccicons}

\fancypagestyle{firststyle}{%
  \fancyhf{}%
  \fancyfoot[L]{\footnotesize %
	\ccby~2022. This version of the article is made available under the CC-BY 4.0 license \url{http://creativecommons.org/licenses/by/4.0/}
	and has been accepted for publication, after peer review but is not the Version of Record and does not reflect post-acceptance improvements, or any corrections. The Version of Record is available online at: \url{http://dx.doi.org/10.1007/s00010-022-00876-4}
	}%
}

\usepackage{hyperref}
\definecolor{darkgreen}{rgb}{0,0.4,0}
\definecolor{BrickRed}{rgb}{0.65,0.08,0}
\hypersetup{colorlinks=true,linkcolor=blue,citecolor=red,filecolor=BrickRed,urlcolor=darkgreen}
\newcommand{\Sc}{\mathcal{S}}

\newcommand{\N}{\mathbb{N}}
\newcommand{\Q}{\mathbb{Q}}
\newcommand{\Z}{{\mathbb Z}}
\newcommand{\R}{{\mathbb R}}

\newtheorem{theo}{Theorem}[section]
\newtheorem{lemma}[theo]{Lemma}

\newtheorem{definition}[theo]{Definition}

\theoremstyle{remark}
\newtheorem{remark}[theo]{Remark}

\newcommand{\gA}{A}
\newcommand{\gB}{B}
\newcommand{\gC}{C}

\begin{document}

\author{
Michael Wallner%
\thanks{Michael Wallner was supported by the Austrian Science Fund (FWF):~P~34142.}
}

\newcommand{\addorcid}[1]{\protect\includegraphics[height=3mm]{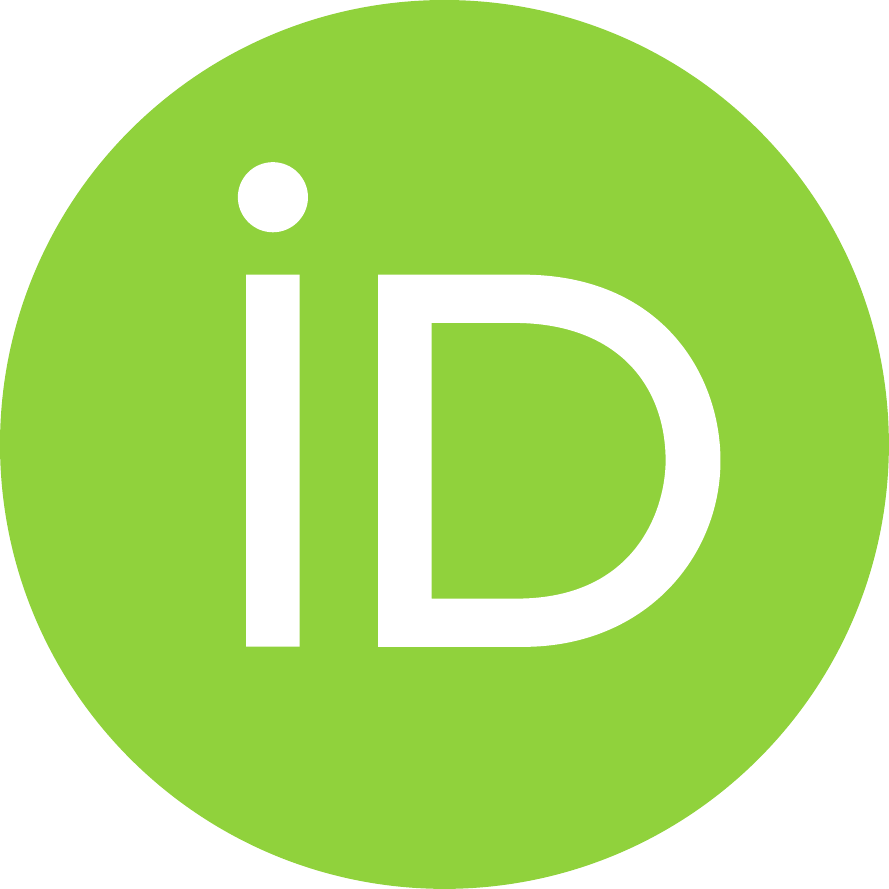} \url{https://orcid.org/#1}}
\newcommand{\Addresses}{{
  \bigskip
  \footnotesize
	
  M.~Wallner, \textsc{%
	 Institute of Discrete Mathematics and Geometry, TU Wien, Wiedner Hauptstra{\ss}e 8--10, 1040 Wien, Austria}\par\nopagebreak
  \textit{E-mail address:} \texttt{michael.wallner@tuwien.ac.at}\par\nopagebreak
  \textit{Website:} \url{https://dmg.tuwien.ac.at/mwallner/} \par\nopagebreak
  \addorcid{0000-0001-8581-449X}

}}

\title{On the critical exponents of generalized ballot sequences in three dimensions and large tandem walks}

\maketitle
\thispagestyle{firststyle}

\renewenvironment{abstract}
{\small
\begin{center}
\bfseries \abstractname\vspace{-.5em}\vspace{0pt}
\end{center}
\list{}{%
	\setlength{\leftmargin}{7mm}
	\setlength{\rightmargin}{\leftmargin}%
}%
\item\relax}
{\endlist}

\begin{abstract}
	We answer some questions on the asymptotics of ballot walks raised in [S.\ B.\ Ekhad and D.\ Zeilberger, April 2021] and prove that these models are not D-finite. 
	This short note demonstrates how the powerful tools developed in the last decades on lattice paths in convex cones help us to answer some challenging problems that were out of reach for a long time. 
	On the way we generalize tandem walks to the family of large tandem walks whose steps are of arbitrary length and map them bijectively to a generalization of ballot walks in three dimensions.
	
	\medskip

\end{abstract}

\section{Introduction}
\label{sec:intro}

In a recent paper~\cite{EkhadZeilberger2021Deep} Ekhad and Zeilberger numerically investigated several problems on the asymptotic enumeration of lattice paths. 
This has led to many conjectures, among which we will answer some in this note. 
They are all related to generalizations of the famous ballot problem which we present now.

The \emph{ballot problem} is probably one of the oldest problems on lattice path enumeration.
There exist various formulations, yet all of them have the following question in common:
Consider an election between two candidates, say A and B, who both receive $n$ votes. 
What is the probability that A is never behind B during the count?
The result is surprisingly simple:
\begin{align*}
	\frac{1}{n+1}.
\end{align*}

The first reported answer was given by Bertrand~\cite{Bertrand87} in 1887, which was shortly after also rigorously proved by Andr\'e~\cite{Andre87}.
Let us call the sequence of possible countings leading to a tie a \emph{ballot sequence}.
Then, what lies behind this beautiful result, is the fact that the total number of votes are counted by the famous \emph{Catalan numbers} 
\begin{align}
	\label{eq:Catalan}
	\frac{1}{n+1} \binom{2n}{n}.
\end{align}

A proof strategy consists in reformulating this question into one about lattice paths in the integer lattice $\Z^2$. 
We model a vote for A by a unit step in $x$-direction $(1,0)$ and a vote for B by a unit step in $y$-direction $(0,1)$. 
The paths start at the origin $(0,0)$ and end at $(n,n)$. 
We are then interested in the total number of paths that are constrained to the region 
\begin{align*}
	\left\{ (x,y) \in \Z^2~:~x \geq y \right\}. 
\end{align*}
These paths are also famously known as \emph{Dyck paths}, being obviously enumerated by the Catalan numbers~\cite{stan99}.
For more on the ballot problem and the history of lattice paths we refer to the beautiful surveys by Humphreys~\cite{Humphreys2010History} and Krattenthaler~\cite{Krattenthaler2015LP}. 

As an answer to such questions, one strives for the most detailed possible one. 
Most desirable is therefore a closed-form expression such as~\eqref{eq:Catalan}.
Often, however, this does not exist or is very complicated but we are still sometimes able to derive its asymptotics. 
Let $a_n$ be the number of paths of length $n$. 
Then for ``well-behaved'' sequences we expect for $n\to\infty$ the structure\footnote{Two sequences $(a_n)_{n \geq 0}$ and $(b_n)_{n \geq 0}$ are asymptotically equivalent, denoted by $a_n \sim b_n$, if and only if $\lim_{n \to \infty} \frac{a_n}{b_n}=1$.}
\begin{align}
	\label{eq:asymptotics}
	a_n \sim C \, \mu^n \, n^\alpha,
\end{align}
where $C$ is a constant, $\mu$ the \emph{exponential growth}, and $\alpha$ the \emph{critical exponent}.
Note that general sequences may include additional terms such as $n!$, $n^n$, or more complicated terms such as $\mu^{n^{\sigma}}$ known as stretched exponentials~\cite{ElveyPriceFangWallner2019Compacted}. 
In the context of this note, the structure~\eqref{eq:asymptotics} will suffice.

Now the ballot problem can be generalized in many ways. 
For our purposes, let $a,b \in \Z_{>0}$ be two positive integers that are relatively prime. 
Then, we consider walks with positive unit steps as above that go from the origin $(0,0)$ to the point $(an,bn)$ and are constrained to the region $bx - ay \geq 0$.  
The special case $a=b=1$ gives the ballot sequence enumerated by Catalan numbers, while $a=1$ and general $b$ leads to \emph{Fuss--Catalan numbers}~\cite{Krattenthaler2015LP}.
The asymptotics of the general problem is answered in~\cite{BanderierWallner16} and can be tricky due to possible cancellations (see for example the conjecture by Knuth\footnote{\url{https://www-cs-faculty.stanford.edu/~knuth/flaj2014.pdf} [viewed 25.05.2021]} for $a=2$ and $b=5$ answered therein).
What is remarkable in these 2-dimensional models is the \emph{universality} of the critical exponent~$\alpha$, in the sense that it is independent of the chosen model and always 
$-3/2$.

\newcommand{\bara}{\gA}
\newcommand{\barb}{\gB}
\newcommand{\barc}{\gC}

A natural next step is to ask the same question in higher dimensions. 
Here the critical exponent is not universal as we will see below and many open questions remain. 
From now on we focus on \emph{generalized ballot sequences in 3 dimensions} which count the following:
\begin{definition}[Generalized 3-ballot walks]
\label{def:gen3ballot}
Let $a,b,c \in \Z_{>0}$ be positive integers 
such that $\gcd(a,b,c)=1$ 
and $M=\operatorname{lcm}(a,b,c)$. 
Generalized 3-ballot walks are walks from the origin $(0,0,0)$ to $(an, bn, cn)$ using unit steps $(1,0,0), (0,1,0), (0,0,1)$ that are constrained to stay in the region
\begin{align*}
	\bara x &\geq 
	\barb y \geq 
	\barc z \geq 0, 
	&& \text{ where } &
	\bara&:=\frac{M}{a},~
	\barb:=\frac{M}{b},~
	\barc:=\frac{M}{c}.
\end{align*}
\end{definition}

Ekhad and Zeilberger~\cite{EkhadZeilberger2021Deep} numerically computed the critical exponent~$\alpha$ for all models with $1 \leq a \leq b \leq c \leq 4$. 
We will derive closed forms for these values, show that all except $a=b=c=1$ are irrational, and explain why this suffices to show that the corresponding generating function is not D-finite.
Note that a generating function is \emph{D-finite} if it satisfies a linear differential equation with polynomial coefficients. This is equivalent to its counting sequence being \emph{P-recursive}, i.e.\ it satisfies a linear recurrence relation with polynomial coefficients. 
For more information on these notions see~\cite{stan99}. 

\begin{remark}
	The list in~\cite[page~5]{EkhadZeilberger2021Deep} unfortunately contains some typos. It does not state $(a,b,c)$ but $(\bara,\barb,\barc)$. Note that it may also be interpreted as $(\barc,\barb,\bara)$; see Lemma~\ref{lem:swapAC}. 
\end{remark}

\section{Closed forms for the numerical estimates of Ekhad--Zeilberger}

\begin{figure}[!ht]
	\centering
	\scalebox{1.1}{
%
%

\begin{tikzpicture}[scale=0.5]
	
	\newcommand{\Asc}{3}
	\newcommand{\Bsc}{2}
	\newcommand{\Csc}{1}

	\newcommand*{\quadw}{14}
	\newcommand*{\quadh}{11}
	\pgfmathsetmacro{\dd}{\quadh-\quadw}
	
	\draw[step=1cm, Gray,very thin,dashed] (0,0) grid (\quadw-0.1,\quadh-0.1);
	
	%
		%
	
	\draw [<->] (0,-2) -- (0,\quadh);
	\draw [<->] (-2,0) -- (\quadw,0);

	\newcommand{\drawcrossing}[2]{
	\coordinate (a) at #1;
	\coordinate (b) at #2;
	\draw[-,white,decorate,double=celestialblue,line width=1.8, double distance=1.8] {($(a)!0.2!(b)$) -- ($(a)!0.8!(b)$)};
	\draw[-,color=celestialblue,line width=1.8] {($(a)!0.12!(b)$) -- ($(a)!0.22!(b)$)};
	\draw[-,color=celestialblue,line width=1.8] {($(a)!0.78!(b)$) -- ($(a)!0.82!(b)$)};
	};
	
	\def\xe{0}
	\def\ye{0}
	\coordinate (current point) at (0,0);
	\foreach \x/\y in {%
		\Asc/0,
		-\Bsc/\Bsc,
		\Asc/0,
		\Asc/0,
		0/-\Csc,
		0/-\Csc
  }
	{	
		\pgfmathsetmacro{\xs}{\xe}
		\pgfmathsetmacro{\ys}{\ye}		
		
		\pgfmathsetmacro{\xe}{\xs+\x}
		\pgfmathsetmacro{\ye}{\ys+\y}
		
		\draw [->,line width=1.8,color=celestialblue] (\xs,\ys) -- (\xe,\ye);
		
		\xdef\xe{\xe}
		\xdef\ye{\ye}
	}	;
	
	\drawcrossing{(6,1)}{(4,3)};
	
	\foreach \x/\y in {%
		-\Bsc/\Bsc,
		-\Bsc/\Bsc,
		-\Bsc/\Bsc,
		0/-\Csc,
		0/-\Csc,
		\Asc/0,
		-\Bsc/\Bsc,
		-\Bsc/\Bsc,
		\Asc/0,
		\Asc/0,
		0/-\Csc,
		0/-\Csc,
		0/-\Csc,
		0/-\Csc,
		-\Bsc/\Bsc
  }
	{	
		\pgfmathsetmacro{\xs}{\xe}
		\pgfmathsetmacro{\ys}{\ye}		
		
		\pgfmathsetmacro{\xe}{\xs+\x}
		\pgfmathsetmacro{\ye}{\ys+\y}
		
		\draw [->,line width=1.8,color=celestialblue] (\xs,\ys) -- (\xe,\ye);
		
		\xdef\xe{\xe}
		\xdef\ye{\ye}
	}	;
	
	\drawcrossing{(2,4)}{(4,4)};
	\drawcrossing{(3,7)}{(1,9)};
	
	\foreach \x/\y in {%
		-\Bsc/\Bsc,
		-\Bsc/\Bsc,
		0/-\Csc
  }
	{	
		\pgfmathsetmacro{\xs}{\xe}
		\pgfmathsetmacro{\ys}{\ye}		
		
		\pgfmathsetmacro{\xe}{\xs+\x}
		\pgfmathsetmacro{\ye}{\ys+\y}
		
		\draw [->,line width=1.8,color=celestialblue] (\xs,\ys) -- (\xe,\ye);
		
		\xdef\xe{\xe}
		\xdef\ye{\ye}
	}	;
	
	\drawcrossing{(0,9)}{(0,7)};
	
	\foreach \x/\y in {%
		0/-\Csc,
		0/-\Csc,
		0/-\Csc,
		0/-\Csc,
		0/-\Csc,		
		0/-\Csc,
		0/-\Csc,
		0/-\Csc
  }
	{	
		\pgfmathsetmacro{\xs}{\xe}
		\pgfmathsetmacro{\ys}{\ye}		
		
		\pgfmathsetmacro{\xe}{\xs+\x}
		\pgfmathsetmacro{\ye}{\ys+\y}
		
		\draw [->,line width=1.8,color=celestialblue] (\xs,\ys) -- (\xe,\ye);
		
		\xdef\xe{\xe}
		\xdef\ye{\ye}
	}	;
		
	\draw [->,line width=1.8,color=celestialblue] (0,1) -- (0,0);	

	\coordinate (cross) at (\quadw-\Asc-1,\quadh-\Bsc-1);
	\foreach \x/\y in { \Asc/0, -\Bsc/\Bsc, 0/-\Csc }
	{	
		\draw[->, line width=2] (cross) -- ($(cross)+(\x,\y)$);
	}
\end{tikzpicture}

%
}
	\caption{%
	3-ballot walks of the model $(a,b,c)=(2,3,6)$ are in bijection with large tandem excursions in the quarter plane with step set $\{(3,0), (-2,-2), (0,-1)\}$; for more details see Lemma~\ref{lem:bijection}.
	}
	\label{fig:largetandem}
\end{figure}
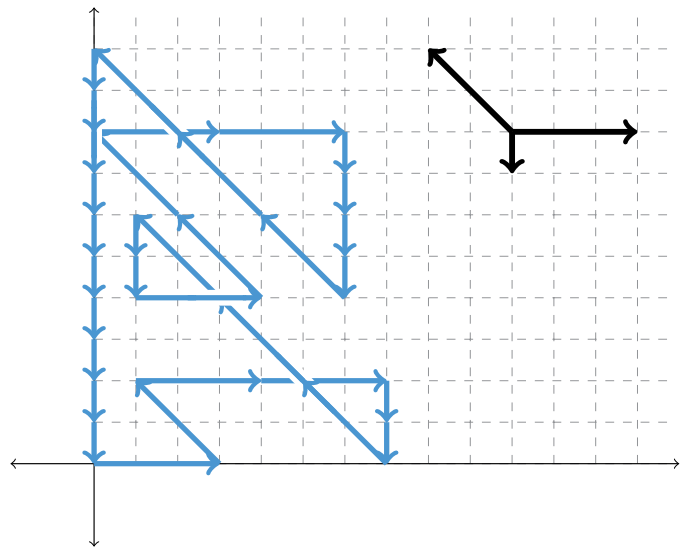

The theory to derive the critical exponent comes from the enumeration of lattice paths in the quarter plane, i.e.\ the first quadrant~\cite{bomi10}. 
As a first step we will map the 3-dimensional walks to the following class of 2-dimensional walks; see Figure~\ref{fig:largetandem}.

\begin{definition}[Large tandem walks]
	\label{def:largetandem}
	Let $\gA,\gB,\gC \in \Z_{>0}$ be positive integers such that $\gcd(\gA,\gB,\gC)=1$. 
\emph{Large tandem walks} are walks starting from the origin $(0,0)$ using the steps $(\gA,0), (-\gB,\gB), (0,-\gC)$ that are constrained to the quarter plane 
\begin{align*}
	\left\{(x,y) \in \Z^2 : x,y\geq 0\right\}.
\end{align*}
\end{definition}

\pagebreak
The 2-dimensional walks with the step set $\{(1,0),(-1,1),(0,-1)\}$ that are constrained to the quarter plane are called \emph{tandem walks} and are connected to a wide range of research problems~\cite{ChyzakYeats2020Tandem,bomi10,biane2021mating}.
The name stems from queueing theory, as they indeed represent the behavior of two $M/M/1$ queues in ``tandem''; see~\cite[Section~4.7]{Fayolle}.

The general step set $\{ (\bara,0)$, $(-\barb,\barb)$, $(0,-\barc) \}$ consists for $(\bara,\barb,\barc)\neq(1,1,1)$ of \emph{large steps}.
Note that the steps $\{-1,0,1\}^2$ are called \emph{small steps}, all others are called large steps.
Walks with large steps are substantially harder to count and the enumeration problem is still widely open. 
A lot of progress has recently been made by Bostan, Bousquet-M\'elou, and Melczer~\cite{Bostan2021Large}, who developed a general theory for walks with large steps in the quarter plane.
We will use this theory to derive closed forms for $\alpha$. 
Remarkably, all that can be done fully automatically, as explained in \cite[Section~8.2.1]{Bostan2021Large}. 

The following result gives an explicit bijection between generalized 3-ballot walks and large tandem excursions, where \emph{excursions} are walks from the origin to the origin.

\begin{lemma}
	\label{lem:bijection}
	The map $\varphi : \Z^3 \to \Z^2$ defined by
	\begin{align*}
		\varphi : (x,y,z) \mapsto \left( \bara x - \barb y, \barb y - \barc z \right)
	\end{align*}
	is a bijection between 
	generalized 3-ballot walks
	from Definition~\ref{def:gen3ballot}
	and large tandem excursions of length $3n$ from Definition~\ref{def:largetandem}.
\end{lemma}

\begin{proof}
	First, note that $M=\operatorname{lcm}(a,b,c)=\operatorname{lcm}(\bara,\barb,\barc)$ due to $\gcd(a,b,c)=\gcd(\bara,\barb,\barc)=1$. 
	Thus, there is a one-to-one correspondence between $(a,b,c)$ and $(\bara,\barb,\barc)$: $a\bara = b\barb = c\barc = M$.
	Second, note that instead of the absolute coordinates, the map keeps track of the weighted relative differences of the coordinates. 
	Thereby, the space constraint $\bara x \geq \barb y \geq \barc z \geq 0$ implies that each 2D-coordinate is nonnegative and therefore lies in the quarter plane. 
	Plugging the unit vectors into $\varphi$ gives the new steps, while the starting point $(0,0,0)$ and the end point $(an,bn,cn)$ both give $(0,0)$. 
	As there is a one-to-one correspondence between the steps this gives the claimed bijection.
\end{proof}

Before we continue we need to introduce some notation.
We associate with a stepset $\Sc \subset \Z^2$ a step polynomial
$
	S(x,y) = \sum_{(i,j) \in \Sc} x^i y^j.
$
This gives in our case 
\begin{align}
	\label{eq:steppoly}
	S(x,y) = x^{\bara} + \left(\frac{y}{x}\right)^{\barb} + \frac{1}{y^{\barc}}.
\end{align}
Let $e_n$ be the number of excursions of length $n$ and let $E(t) = \sum_{n \geq 0} e_n t^n$ be the generating function of excursions. 
%
We define the period $p$ of excursions as $ p := \gcd\{n \in \N : e_n \neq 0\}$.
In other words, it holds that $E(t) = \sum_{m \geq 0} e_{pm} t^{pm}$.
Observe that for large tandem excursions we have $p = a+b+c$. 
Moreover, the origin $(0,0)$ is \emph{reachable from infinity} if there exists a quadrant walk that starts from a point $(k,\ell) \in \Z_{>0}^2$ and ends at $(0,0)$.
Note that for large tandem excursions this is the case, e.g., for the point $(\gB,\gB\gC-\gB)$.

The main result we will need is the following adapted theorem which we state for excursions only.
%
%
%
\begin{theo}[Adapted {\cite[Theorem~7 and Theorem~11]{Bostan2021Large}}; see also~\cite{DenisovWachtel15}]
	\label{theo:mainbostanetal}
	Let $\Sc \subset \Z^2$ be a step set that is not contained in a half-plane and contains an element of $\N^2$. Then the step polynomial $S(x,y)$ has a unique critical point $(X,Y)$ in $\R_{>0}^2$ (that is, a solution of $S_x(X,Y) = S_y(X,Y) = 0$), which satisfies $S_{xx}(X,Y) > 0$ and $S_{yy}(X,Y) >0$. Define
	\begin{align}
		\notag
		\mu &= S(X,Y), \\
		\label{eq:c}
		\gamma &= \frac{S_{xy}(X,Y)}{\sqrt{S_{xx}(X,Y)S_{yy}(X,Y)}}, \\
		\notag
		\alpha &= -1 - \frac{\pi}{\arccos(-\gamma)}.
	\end{align}
	Furthermore, let $(0,0)$ be reachable from infinity. 
	If the period of excursions satisfies $p=1$, then there exists a positive constant $\kappa$ such that, as $n$ goes to infinity
	\begin{align*}
		e_{n} &\sim \kappa \, \mu^{n} \, {n}^{\alpha}.
	\end{align*}
	If the period of excursions satisfies $p>1$, then we have for $n=pm$ large enough\footnote{For two sequences $(a_n)_{n \geq 0}$ and $(b_n)_{n \geq 0}$ we write $a_n = \Theta(b_n)$, if and only if there exist constants $c_1, c_2 >0$ and $N \in \N$ such that $c_1 |b_n| \leq |a_n| \leq c_2 |b_n|$ for $n \geq N$.}
	\begin{align*}
		e_{n} &= \Theta\left( \mu^{n} {n}^{\alpha} \right).
	\end{align*}
	Then $\alpha$ is called the \emph{critical excursion exponent}.
	If $\alpha \notin \Q$ (or, equivalently, $\arccos(-\gamma)$ is not a rational multiple of $\pi$), then 
	the series $E(t)$ is not D-finite.
\end{theo}

\begin{remark}
	In~\cite{Bostan2021Large} the authors consider the generating function $Q(x,y;t) = \sum_{i,j,n \geq 0} q_{ijn} x^i y^j t^n$ where $q_{ijn}$ is equal to the number of walks from $(0,0)$ to $(i,j)$ of length $n$ restricted to the quarter plane. Then we have $E(t) = Q(0,0;t)$ and as D-finite functions are closed under algebraic substitutions (see, e.g., \cite{flaj09}), $Q(x,y;t)$ is not D-finite either if $E(t)$ is not D-finite.
\end{remark}

The authors of~\cite{Bostan2021Large} have implemented the procedure to determine $\gamma$ and $\alpha$ in a publicly available Maple worksheet\footnote{\url{https://www.labri.fr/perso/bousquet/publis.html} [viewed 25.05.2021]}.
More details of the involved methods are explained in \cite[Section~8.2.1]{Bostan2021Large}.
Using the associated $\texttt{minPs}$ command we compute the closed forms of the critical exponents conjectured in~\cite{EkhadZeilberger2021Deep} and list them in Table~\ref{tab:alpha}.
Now using the arguments from~\cite{BostanRaschelSalvy2014Nondfinite} (which are also used in~\cite{Bostan2021Large}) we directly get that the critical exponent $\alpha$ is irrational in all cases, except for $(a,b,c)=(1,1,1)$, and hence the generating function $E(t)$ (as well as $E(t^{1/p})$) is not D-finite.
For the case $(a,b,c)=(1,1,1)$ of tandem excursions it is known that the generating function is D-finite (even hypergeometric); see~\cite[Proposition~9]{bomi10} and \cite{GesselZeilberger92Weyl}.

\begin{table}[!ht]
\renewcommand{\arraystretch}{1.15}
\centering
\begin{tabular}{ccccc}
\toprule
\multicolumn{2}{c}{Model} & \multicolumn{3}{c}{$\alpha$} 
\\\cmidrule(lr){1-2} \cmidrule(lr){3-5}
$(a,b,c)$ & 
$(\bara,\barb,\barc)$ &
Conjecture~\cite{EkhadZeilberger2021Deep} & 
Correct Approx. &
Closed form \\
\midrule
$(1,1,1)$ & 
$(1,1,1)$ & 
$-4$ &
$-4$ &
$-4$ \\
$(1,2,2)$ &
$(2,1,1)$ &
$-3.7312$ &
$-3.7312$ &
$-1 - \frac{\pi}{\arccos(1/\sqrt{6})}$\\
$(1,1,2)$ &
$(2,2,1)$ &
$-4.2884$ &
$-4.28854$ &
$-1 - \frac{\pi}{\arccos(1/\sqrt{3})}$\\
$(1,3,3)$ &
$(3,1,1)$ &
$-3.5976$ &
$-3.59758$ &
$-1 - \frac{\pi}{\arccos(1/\sqrt{8})}$\\
%
$(2,3,6)$ & $(3,2,1)$ & $-4.055$ & $-4.05556$ & $-1-\frac{\pi}{\arccos(2/\sqrt{15})}$ \\ 
$(2,3,3)$ & $(3,2,2)$ & $-3.8375$ & $-3.83755$ & $-1-\frac{\pi}{\arccos(1/\sqrt{5})}$ \\ 
$(1,1,3)$ & $(3,3,1)$ & $-4.4455$ & $-4.44572$ & $-1-\frac{\pi}{\arccos(\sqrt{3/8})}$ \\ 
$(2,2,3)$ & $(3,3,2)$ & $-4.1695$ & $-4.16962$ & $-1-\frac{\pi}{\arccos(\sqrt{3/10})}$ \\ 
$(1,4,4)$ & $(4,1,1)$ & $-3.515$ & $-3.51519$ & $-1-\frac{\pi}{\arccos(1/\sqrt{10})}$ \\ 
$(1,2,4)$ & $(4,2,1)$ & $-3.9091$ & $-3.90911$ & $-1-\frac{\pi}{\arccos(\sqrt{2/9})}$ \\ 
$(3,4,12)$ & $(4,3,1)$ & $-4.2454$ & $-4.24544$ & $-1-\frac{\pi}{\arccos(3/\sqrt{28})}$ \\ 
$(3,4,6)$ & $(4,3,2)$ & $-4.0237$ & $-4.02370$ & $-1-\frac{\pi}{\arccos(3/\sqrt{35})}$ \\ 
$(3,4,4)$ & $(4,3,3)$ & $-3.8834$ & $-3.88346$ & $-1-\frac{\pi}{\arccos(\sqrt{3/14})}$ \\ 
$(1,1,4)$ & $(4,4,1)$ & $-4.5453$ & $-4.54551$ & $-1-\frac{\pi}{\arccos(\sqrt{2/5})}$ \\ 
$(3,3,4)$ & $(4,4,3)$ & $-4.12019$ & $-4.12021$ & $-1-\frac{\pi}{\arccos(\sqrt{2/7})}$ \\

\bottomrule
\end{tabular}
\vspace{-1mm}
\caption{Proved closed forms of the critical exponent $\alpha$ of excursions for the conjectured values from~\cite[page~5]{EkhadZeilberger2021Deep}. 
As always, $A=M/a$, $B=M/b$, $C=M/c$, and $M=\operatorname{lcm}(a,b,c)$.
Note that $\alpha$ is irrational in all cases except the first, and therefore the excursions' generating function is \emph{not} D-finite.}
\label{tab:alpha}
\end{table}

This gives the closed-form expressions requested in~\cite{EkhadZeilberger2021Deep} and thereby makes the OEIS USD~100 richer. 
In the next section we go one step further and investigate the models associated with generalized 3-ballot walks, namely large tandem excursions, and uncover a surprisingly simple structure for~$\gamma$. 
This gives a tool to prove the non-D-finiteness of many such models.

\vspace{-2mm}
\section{Generalized 3-ballot walks and large tandem walks}
\vspace{-1mm}

%
%

As a first observation we note that the number of excursions is the same in the models $(\bara,\barb,\barc)$ and $(\barc,\barb,\bara)$. This follows from the following lemma.

\begin{lemma}
	\label{lem:swapAC}
	There exists a bijection between large tandem excursions in the models $(\bara,\barb,\barc)$ and $(\barc,\barb,\bara)$.
\end{lemma}

\begin{proof}
	Let a large tandem excursion in the model $(\bara,\barb,\barc)$ be given. 
	It is built from the step set $\{(\gA,0), (-\gB,\gB), (0,-\gC)\}$.
	First, we reverse the time by reading the steps backwards. This gives the new step set $\{ (-\gA,0)$, $(\gB,-\gB), (0,\gC) \}$.
	Second, we rotate along $y=x$, and get the step set $\{ (\gC,0)$, $(-\gB,\gB), (0,-\gA) \}$ corresponding to the model $(\barc,\barb,\bara)$. 
	Note that both operations transform excursions into excursions and remain in the quarter plane. 
	For more details and further transformations see~\cite{BanderierLacknerWallner2020Latticepathology}.
\end{proof}

Next we compute the critical excursion coefficient using Theorem~\ref{theo:mainbostanetal}. 
Observe the invariance in swapping $\bara$ and $\barc$ implied by Lemma~\ref{lem:swapAC}.

\begin{lemma}
	\label{lem:clargetandem}
	The exponential growth $\mu$, the constant $\gamma$, and the critical exponent $\alpha$ from Equation~\eqref{eq:c} for large tandem excursions with the step set $\{ (\gA,0), (-\gB,\gB), (0,-\gC) \}$ are equal to
\begin{align}
	\mu &= \gC \left(\frac{\gA^\gB \gB^\gA}{\gC^{\gA+\gB}}\right)^{\gC/(\gA\gB+\gA\gC+\gB\gC)} \left(\frac{1}{\gA}+\frac{1}{\gB}+\frac{1}{\gC}\right), \notag\\
	\gamma &= -\frac{\gB}{\sqrt{(\gA+\gB)(\gB+\gC)}}, \label{eq:cgeneric} \\
	\alpha &= -1 - \frac{\pi}{\arccos\left( \gB/\sqrt{(\gA+\gB)(\gB+\gC)} \right)}. \notag
\end{align}
\end{lemma}

\begin{proof}
	The step polynomial of large tandem walks is given in~\eqref{eq:steppoly}:
	$
		S(x,y) = x^{\gA} + (\frac{y}{x})^{\gB} + \frac{1}{y^{\gC}}.
	$
	According to Theorem~\ref{theo:mainbostanetal} we need to compute the unique values $X,Y>0$ such that $S_x(X,Y)=S_y(X,Y)=0$. This gives the following equalities
	\begin{align*}
		\gA X^{\gA} = \gB \left(\frac{Y}{X}\right)^{\gB} = \gC \frac{1}{Y^{\gC}}, \qquad \text{ where }
	\end{align*}	
	\begin{align*}
		X &= \left(\frac{\gB^\gC \gC^\gB}{\gA^{\gB+\gC}}\right)^{1/(\gA\gB+\gA\gC+\gB\gC)} && \text{and} &
		Y &= \left(\frac{\gC^{\gA+\gB}}{\gA^\gB \gB^\gA}\right)^{1/(\gA\gB+\gA\gC+\gB\gC)}.
	\end{align*}	
	Using these relations we compute and simplify $S_{xy}(X,Y)$, $S_{xx}(X,Y)$, and $S_{yy}(X,Y)$:
	\begin{align*}
		S_{xy}(X,Y) &= -\gB^2 \frac{Y^{\gB-1}}{X^{\gB+1}}, \\
		S_{xx}(X,Y) &= \gA(\gA-1)X^{\gA-2} + \gB(\gB+1)\frac{Y^{\gB}}{X^{\gB+2}} 
		             = (\gA+\gB) \gB \frac{Y^{\gB}}{X^{\gB+2}}, \\
		S_{yy}(X,Y) &= \gB(\gB-1)\frac{Y^{\gB-2}}{X^{\gB}} + \gC(\gC+1)\frac{1}{Y^{\gC+2}} 
		             = (\gB+\gC) \gB \frac{Y^{\gB-2}}{X^{\gB}}.
	\end{align*}
	Thus, using Formula~\eqref{eq:c} we get the claimed closed forms. Note that the specific values of $X$ and $Y$ are only needed for the determination of $\mu$.
\end{proof}

This remarkably simple formula allows us to determine many non-D-finite models.

\begin{theo}
	\label{theo:main}
	The generating function of large tandem walks (and in particular large tandem excursions) with parameters $(\gA,\gB,\gC)$ is not D-finite if 
	\[
		\frac{\gB^2}{(\gA+\gB)(\gB+\gC)} \notin \left\{ \frac{1}{4}, \frac{1}{2}, \frac{3}{4}\right\}.
	\]
\end{theo}

\begin{proof}
Building on the results of Lemma~\ref{lem:clargetandem}, 
it remains to determine whether $\alpha \in \Q$, 
which is equivalent to $\arccos(-\gamma)/\pi \in \Q$.
There exist general classification results for the case that $\gamma$ is a square root of a rational~\cite{Varona2006arccosine}: For $r \in \Q$ with $0 < r < 1$, the value $\arccos(\sqrt{r})/\pi$ is rational if and only if $r \in \{1/4, 1/2, 3/4\}$. Hence, in our case everything depends on $r = \gamma^2 = \frac{\gB^2}{(\gA+\gB)(\gB+\gC)}$. 
\end{proof}

All possible rational values of $\alpha$ are summarized in Table~\ref{tab:3ballotalphagen}.
It shows the only possible candidates for D-finite generating functions; all other cases are therefore proved to be non-D-finite.

\begin{table}[!ht]
\renewcommand{\arraystretch}{1.25}
\begin{center}
\begin{tabular}{cccp{0.2cm}c}
\toprule
$\gamma^2$ & $\arccos(-\gamma)$ & $\alpha$ && Large tandem walks $(\gA,\gB,\gC)$
\\
\midrule
$1/4$ & $\pi/3$ & $-4$ && 
	\begin{tabular}{c}
		$(1,1,1)$, $(1,3,6)$, $(2, 14, 35)$, 
		$(3, 6, 10)$, $(3, 15, 35)$,  \dots
	\end{tabular} \\
%
%
$1/2$ & $\pi/4$ & $-5$ && 
	\begin{tabular}{c}
		$(2, 6, 3)$, $(3, 15, 10)$, $(4, 28, 21)$, 
		$(5, 20, 12)$, $(5, 45, 36)$, \dots 
	\end{tabular} \\
%
%
$3/4$ & $\pi/6$ & $-7$ && 
	\begin{tabular}{c}
		$(4, 60, 15)$, $(5, 45, 9)$, $(6, 42, 7)$, 
		$(7, 189, 54)$, $(9, 99, 22)$, \dots
	\end{tabular} 
%
%
\\
\bottomrule
\end{tabular}
\caption{By Theorem~\ref{theo:main} large tandem excursions and the bijectively equivalent generalized $3$-ballot walks are not D-finite if $\gamma^2 = \frac{\gB^2}{(\gA+\gB)(\gB+\gC)} \notin \left\{ \frac{1}{4}, \frac{1}{2}, \frac{3}{4}\right\}$. The table lists the only $3$ exceptional values for $\gamma^2$, the associated rational critical excursion exponents $\alpha$, and examples exhibiting this behavior.
The generating function of classical tandem excursions $(1,1,1)$ is D-finite, while the nature of the other models is unknown.}
\label{tab:3ballotalphagen}
\end{center}
\end{table}



The question remains if the models associated with rational $\alpha$ are D-finite or not. 
Observe that they cannot be algebraic as $\alpha$ is a negative integer~\cite{flaj09}.
Furthermore, note that the non-D-finiteness results from~\cite{Bostan2021Large} simultaneously show that the orbits of the models are infinite. 
The equivalence of these two statements was proved in the case of models with small steps, yet it is conjectured to hold in greater generality.
The authors provide the command \texttt{evalPQ} to possibly find a lower bound on the orbit size as well as the command \texttt{PQList} to compute the orbit elements.
Now, using these commands we find that the classical tandem model $(1,1,1)$ (whose generating function is of course D-finite) has a finite orbit of size~$6$, while for all other models from Table~\ref{tab:3ballotalphagen} we could not find any lower bounds on the orbit size as the \texttt{evalPQ} aborted after a certain threshold, which we set to $1000$ and which is related to the orbit size. 
Thus, the orbits could be infinite or just very large.
As mentioned above, a finite orbit does not directly prove D-finiteness, 
however, in the case of models with small steps it was proved that the model is D-finite if and only if the orbit is finite. 

Let us end with some simple observations regarding the models with rational $\alpha$.
First, it is easy to see that fixing two of the three parameters $\gA,\gB,\gC$ completely determines the third one.
Second, for $\gamma^2=1/4$ it holds
 that $\gA<\gB<\gC$ or $\gC<\gB<\gA$, while for the other two cases $\gamma^2=1/2$ and $\gamma^2=3/4$ we have $\gA,\gC < \gB$.
Third, these families are infinite. 
In the first case, consider for example the parameters $\left(\gA,(\gA-1)\gA,(\gA-1)(3\gA-4)\right)$ whose greatest common divisor is equal to one for odd values of $\gA>1$. 
In the second case consider $\left(\gA,(\gA-1)\gA,(\gA-1)(\gA-2)\right)$ for odd $\gA>1$, and in the third case $\left(\gA,(\gA-1)\gA,\frac{(\gA-1)(\gA-4)}{3}\right)$ for $\gA=6k+1$ with $k>0$.

\section{Future work}
\label{sec:future}

We conclude this short note with some interesting questions for future research. 
Common to all is the need to develop general methods for the extremely useful lattice paths (with large steps) in the quarter plane and higher dimensions.

Firstly, it would be interesting to investigate the simplest model of each case in Table~\ref{tab:3ballotalphagen}.
We have not been able to do this so far, and it seems that we have reached the state-of-the-art of our current understanding of such lattice paths. 
Moreover, large tandem walks by themselves seem to be quite interesting and deserve a further study. 
E.g., what is the nature/asymptotics of $Q(1,1;t)$? 
For classical tandem walks $Q(1,1;t)$ is even algebraic, while $Q(0,0,t)$ is D-finite transcendental~\cite[Proposition~9]{bomi10}. However, we have not been able to guess any algebraic or differential equation for $Q(1,1;t)$ using the first $300$ terms for the other models in Table~\ref{tab:3ballotalphagen} and any model satisfying $1 \leq \gA,\gB,\gC \leq 10$ with $\gcd(\gA,\gB,\gC)=1$.

Secondly, it remains to consider higher dimensional $d$-ballot paths, generalizing Definition~\ref{def:gen3ballot}. 
The bijection of Lemma~\ref{lem:bijection} directly generalizes and leads to excursions in $d-1$~dimensions that are constrained to stay in the positive orthant.
Zeilberger~\cite{EkhadZeilberger2021Deep} experimentally considered several models in $4$~dimensions and gave more conjectures for their asymptotics. 
However, our results rely heavily on the deep connections summarized in Theorem~\ref{theo:mainbostanetal}, and we are lacking such an understanding of higher dimensional paths at the moment. 
Note however that for 3-dimensional walks constrained to an octant it may be possible to use the ideas from~\cite{Bogosel2020}, in which some results on the associated critical exponents are available.

Thirdly, classical tandem walks are in bijection with many other combinatorial models, such as Young tableaux with at most $3$ rows or Motzkin paths, to name a few. 
Do these bijections generalize to large tandem walks, or are there similar bijections?

%

\subsection*{Acknowledgement}

I would like to thank Mireille Bousquet-Mélou for useful comments on the periodic non-D-finiteness results in~\cite{Bostan2021Large}.
I also would like to thank the anonymous referee for his or her great interest in this work and the detailed comments and questions that inspired much of Section~\ref{sec:future}.

\addcontentsline{toc}{section}{References}
\bibliographystyle{cyrbiburl}
\bibliography{mybib_Zeilberger}
\label{sec:biblio}

\Addresses

\end{document}